\documentclass[a4paper,11pt]{article}
\usepackage{amsmath, amsfonts}
\usepackage[amsmath,thmmarks]{ntheorem}

\title{On Sidki's presentation for orthogonal groups}
\author{Justin M\textsuperscript{c}Inroy \\ \footnotesize{Heilbronn Institute for Mathematical Research, School of Mathematics,} \\ \footnotesize{University of Bristol, University Walk, Bristol, BS8 1TW, UK} \\ \footnotesize{Tel: +44 117 954 5661, Fax: +44 117 331 5264,} \\ \footnotesize{email: justin.mcinroy@bristol.ac.uk}
\and Sergey Shpectorov \\ \footnotesize{School of Mathematics, University of Birmingham, Edgbaston, Birmingham,} \\ \footnotesize{B15 2TT, UK, Tel: +44 121 414 6604, Fax: +44 121 414 3389,} \\ \footnotesize{email: S.Shpectorov@bham.ac.uk}}

\newcommand\la{\langle}
\newcommand\ra{\rangle}

\newtheorem{thm}{Theorem}[section]
\newtheorem{lem}[thm]{Lemma}
\newtheorem{cor}[thm]{Corollary}
\newtheorem{prop}[thm]{Proposition}

\theorembodyfont{\normalfont}

\theoremstyle{nonumberplain}
\theoremheaderfont{\normalfont\itshape}
\theorembodyfont{\normalfont}
\theoremseparator{.}
\theoremsymbol{\ensuremath{\square}}
\newtheorem{pf}{Proof}
\qedsymbol{\ensuremath{\square}}
\usepackage[pdftex]{graphicx}

\begin{document}

\maketitle

\begin{abstract}
We study presentations, defined by Sidki, resulting in groups $y(m,n)$ that are conjectured to be finite orthogonal groups of dimension $m+1$ in characteristic two.  This conjecture, if true, shows an interesting pattern, possibly connected with Bott periodicity.  It would also give new presentations for a large family of finite orthogonal groups in characteristic two, with no generator having the same order as the cyclic group of the field.

We generalise the presentation to an infinite version $y(m)$ and explicitly relate this to previous work done by Sidki.  The original groups $y(m,n)$ can be found as quotients over congruence subgroups of $y(m)$.  We give two representations of our group $y(m)$.  One into an orthogonal group of dimension $m+1$ and the other, using Clifford algebras, into the corresponding pin group, both defined over a ring in characteristic two.  Hence, this gives two different actions of the group.  Sidki's homomorphism into $SL_{2^{m-2}}(R)$ is recovered and extended as an action on a submodule of the Clifford algebra.
\end{abstract}

\section{Introduction}
The following is a well-known presentation of the alternating group $A_{m+2}$ which was given by Carmichael in 1923:
\[
A_{m+2} = \la a_1, \dots ,a_m \mid a_i^3=1, (a_ia_j)^2=1,\mbox{ for all } i \neq j \ra.
\]
In 1982, Sidki generalised this to the following:
\[
Y(m,n):= \la a_1, \dots ,a_m \mid a_i^n=1, (a_i^k a_j^k)^2=1,\mbox{ for all } i \neq j, 1 \le k \le n-1 \ra
\]
and he conjectured that these too are all finite groups.

It is clear that $Y(m,3)$ is just Carmichael's presentation for the alternating group.  In \cite{si1, si2}, Sidki identifies the groups $Y(2,n)$, $Y(3,n)$ (when $n$ is odd), $Y(m,2)$ and $Y(m,4)$ and hence shows that they are finite.  In general, however, the question of finiteness remains open.

As well as identifying the groups when either $m$, or $n$ is small, Sidki gives some general results.  He shows that $Y(m,n)$ is perfect provided $m>2$ and $n$ is odd, and that if $n| n'$, then $Y(m,n)$ is a quotient group of $Y(m,n')$ \cite[Theorem A]{si1}.  The second of these results allows a reduction to the cases where $n$ is a prime power.  When $n=2^r$, it is conjectured that $Y(m,2^r)$ is a $2$-group.  If $n$ is odd, however, Sidki shows that $Y(m,n)$ is isomorphic to the following group:
\begin{align*}
y(m,n):= \la a,S_m \mid \, &  a^n=1, [s_1,s_1^{a^k}]=1, 1 \le k \le n-1 \\
& s_1^{1+a+ \dots + a^{n-1}}=1 \mbox{ and }a^{s_i}=a^{-1},2\le i \le m-1 \ra,
\end{align*}
where $S_m$ denotes the symmetric group and $s_i = (i, i+1)$.  

In addition to the results above, a few small cases were resolved by J. Neubuser, W. Felsch and E. O'Brian by direct calculations using the Todd-Coxeter algorithm (see Table \ref{table1}).  These calculation suggest that, for odd $n$, $y(m,n)$ are orthogonal groups in characteristic two of dimension $m+1$ for some suitable quadratic form.  Note that when there is a normal subgroup, namely when $m=2 \mod 4$, this agrees with our assessment and just indicates that the form has a non-trivial radical.

\begin{table}[ht]
\begin{center}
\begin{tabular}{ll}
$Y(3,5) \cong SL_2(16) \cong \Omega^-(4,4)$ & $Y(3,7) \cong \Omega^+(4,8)$\\
$Y(4,5) \cong \Omega(5,4)$& $Y(4,7) \cong \Omega(5,8)$\\
$Y(5,5) \cong \Omega^-(6,4)$& $Y(5,7) \cong \Omega^+(6,8)$\\
$Y(6,5) \cong 4^6:\Omega^-(6,4)$& $Y(6,7) \cong 8^6:\Omega^+(6,8)$\\
$Y(7,5) \cong \Omega^-(8,4)$& \\
$Y(8,5) \cong \Omega(9,4)$& $Y(3,11) \cong \Omega^-(4,32)$\\
$Y(9,5) \cong \Omega^-(10,4)$& $Y(4,11) \cong \Omega(5,32)$\\
$Y(10,5) \cong 4^{10}:\Omega^-(10,4)$& $Y(5,11) \cong \Omega^-(6,32)$
\end{tabular}
\end{center}
\caption{Some small cases}\label{table1}
\end{table}

The picture which emerges from Table \ref{table1} and other known results is quite pretty and somewhat unexpected.  Note that the sequence of groups in Table \ref{table1} indicates a cycle of length four.  Sidki himself conjectured that this is connected with Bott periodicity.  If these presentations were for orthogonal groups, the fact that they do not contain a generator which has the same order as the cyclic group of the field is also unique. 

Out of all the results so far, the most interesting and non-trivial is Sidki's solution of the case $m=3$ and $n$ odd.  He starts by defining the infinite analogue of the group $Y(3,n)$ by letting the $a_i$ have infinite order \cite{si2}: 
$$
Y(3):= \la a_1, a_2, a_3 \mid (a_i^ka_j^k)^2=1,i \neq j, k \in \mathbb{Z}\ra.
$$
He proceeds by identifying $Y(3)$ as a certain subgroup of $SL_2(\mathbb{F}_2[s,s^{-1}])$, where $s$ is an indeterminate.  After doing this, using congruence subgroups, he recovers $Y(3,n)$ as a quotient, hence identifying it as $SL_2(\mathcal{A}_n(\mathbb{F}_2))$, where $\mathcal{A}_n(\mathbb{F}_2)$ is the augmentation ideal in the group algebra $\mathbb{F}_2C_n$.  Note that, when $n$ is odd, $\mathcal{A}_n(\mathbb{F}_2)$ is a direct summand of $\mathbb{F}_2C_n$ and hence is in its own right a commutative ring with one.  Thus, the notation $SL_2(\mathcal{A}_n(\mathbb{F}_2))$ makes sense.  (Recall that $GL_n(R)$ is the group of $n \times n$ matrices with entries in $R$ which are invertible.  If $R$ is commutative, then the determinant can be defined and if further $R$ contains a 1, then $SL_n(R)$ is defined as the kernel of the determinant map.  See, for instance, \cite{omeara} for details.)  In fact, as the augmentation ideal splits as the direct sum of fields, if $n$ is an odd composite, the group over the ring $\mathcal{A}_n(\mathbb{F}_2)$ can be a direct sum of several orthogonal groups over some fields.

The fact that $Y(3)$ is isomorphic to a subgroup of $SL_2(\mathbb{F}_2[s,s^{-1}])$ means that there is a homomorphism $\iota$ from $Y(3)$ to $SL_2(\mathbb{F}_2[s,s^{-1}])$ and this homomorphism is injective.  In this respect, it is worth mentioning that in \cite{si1} Sidki gives a homomorphism from $y(m,n)$ ($\cong Y(m,n)$ when $n$ is odd) to $SL_{2^{m-2}}(\mathbb{F})$, where $\mathbb{F}$ is any field of characteristic two containing the $n$\textsuperscript{th} root of unity.  Clearly $\mathbb{F}$ can be chosen to be a finite field and Sidki's conjecture would follow if only we could show that his homomorphism is injective.

Note the discrepancy between the dimension here, $2^{m-2}$, and the dimension $m+1$ of the groups in Table \ref{table1}.  This suggests that Sidki's action could be like the spin action for some orthogonal group.
 
 In this paper, we work with the infinite version of $y(m,n)$ rather than of $Y(m,n)$ as Sidki does.  We define
$$
y(m):= \la a,S_m \mid [s_1,s_1^{a^k}]=1, k \in \mathbb{Z}, a^{s_i}=a^{-1},2\le i \le m-1 \ra,
$$
where $S_m$ is a symmetric group with generators $s_i =(i,i+1)$ as before.  (In fact, we work with a slightly larger group $\tilde{y}(m) = y(m)\la \tau \ra$, where $\tau$ is the automorphism of $y(m)$ which centralises $S_m$ and inverts $a$.  Sidki describes a similar algebra homomorphism in \cite{si1}.)  The group $y(m)$ is not isomorphic to $Y(m)$; in fact we will show that $y(m)$ is the semidirect product of $Y(m)$ with $S_m$, where $Y(m)$ is defined in the obvious way.  The standard generators $s_i$ of $S_m$ act on $Y(m)$ by permuting and inverting $a_i$ and $a_{i+1}$ and inverting all other $a_j$.  In particular, $\tilde{s}_i = \tau s_i$ just permute the $a_i$ in the natural way.

Sidki's homomorphism for $y(m,n)$ readily generalises to $y(m)$ if we replace the $n$\textsuperscript{th} root of unity with an indeterminate.  This gives a homomorphism $\eta:y(m) \to SL_{2^{m-2}}(\mathbb{F}_2[s,s^{-1}])$.  Naturally, the conjecture is that this homomorphism is injective and so $y(m)$ is isomorphic to the image of $\eta$.  In fact, in the case of $m=3$, the map $\eta$ restricted to $Y(3) < y(3)$ with the choice of generators as given above is the same as the injective homomorphism $\iota$ mentioned before.

Our motivation for this paper comes from Table \ref{table1}.  Instead of $\eta$, we construct, in a very natural and simple way, a homomorphism $\phi$ from the larger group $\tilde{y}(m) := y(m)\la \tau \ra$ to the orthogonal group $SO(V,q)$, where $V$ is a free module of rank $m+1$ over $\mathbb{F}_2[t,t^{-1}]$ and $q$ is a certain quadratic form on $V$.  In fact, it maps into a subgroup generated by orthogonal transvections, which we call $TO(V,q)$.  When $\phi$ is restricted to $y(m)$ it maps into $O'(V,q)$, the group generated by products of an even number of orthogonal transvections.  Note that, although the values of our form $q$ do involve $t^{-1}$, all the matrices in the image of $\phi$ have entries which are polynomials in $t$, hence $\mathrm{Im}(\phi)$ lies in the orthogonal group over $\mathbb{F}_2[t]$.  Our homomorphism $\phi$ very nicely explains the known entries in Table \ref{table1}, as these groups are the homomorphic images of the larger orthogonal group $y(m)$ under the suitable evaluation mappings from $\mathbb{F}_2[t,t^{-1}]$ to the finite fields which send $t$ to the appropriate $n$\textsuperscript{th} of unity.  We hope that this fact and the lower dimension of the representation will make this homomorphism easier to work with and will lead to a solution to Sidki's conjecture.


We construct another homomorphism to the pin group $Pin(V,q)$ for the same $V$ and $q$.  Since we realise the pin group as a group of units in the Clifford algebra $Cl(V,q)$ defined by $V$ and $q$, this endows $\psi(\tilde{y}(m))$ with an action on the Clifford algebra in addition to that of $\phi(\tilde{y}(m))$ on the natural module $V$.  Although the homomorphism $\psi$ adds a layer of abstraction, the definition remains quite natural and easy to define in terms of elements of $Cl(V,q)$.  We note at this point that we need to extend the ring by adding $s:=\sqrt{t}$.

A priori, $\psi$ could be a non-trivial lifting of $\phi$.  Indeed, we show that $Pin(V,q)$ is a non-trivial cover of the subgroup of $TO(V,q) \leq SO(v,q)$ generated by orthogonal transvections if and only if $m$ is even.  However, we show that for our group, $\psi$ is a trivial lifting of $\phi$.  So, $\phi$ and $\psi$ carry the same information.  We conjecture that both $\phi$ and $\psi$ are injective and hence that $\tilde{y}(m)$ is isomorphic to its image under both $\phi$ and $\psi$.

To support this conjecture, we demonstrate the equivalence of our homomorphism $\psi$ to Sidki's $\eta$ in the following way.  After extending the ring further by adding a root $\alpha$ of the minimum polynomial of $a$, we find a submodule of rank $2^{m-2}$ in the Clifford algebra on which the group $\psi(y(m))$ acts, and a particular basis in this submodule so that the corresponding action matrices match Sidki's.  Moreover, we find a submodule of rank $2^{m-1}$ which is invariant under the action of the larger group $\psi(\tilde{y}(m))$.

Furthermore, this opens up the connections with the geometry, namely, with the orthogonal group over $\mathbb{F}_2[s]$ acting on the twin building.  We hope to find in that twin building a simple connected geometry on which $\mathrm{Im}(\psi)$ acts flag-transitively with the amalgam of maximal parabolics encoded in the presentation of $\tilde{y}(m)$.  If this were true, then the injectivity of $\psi$ would follow from Tits' lemma.

In Section \ref{sec:y(m)}, we discuss the relationship between the two infinite groups $y(m)$ and $Y(m)$.  (Recall that the finite groups $y(m,n)$ and $Y(m,n)$ are isomorphic when $n$ is odd.)  We show that $\tilde{y}(m) \cong Y(m) : (S_m \times C_2)$, where $S_m$ permutes naturally the generators $a_1, \dots , a_m$ of $Y(m)$ and the $C_2$ inverts all the generators.  We describe the quadratic form and two homomorphisms in Section \ref{sec:homo}, which is split into three parts.  We begin with some details about our specific form $q$ and then the first subsection has the homomorphism $\phi$ from $\tilde{y}(m)$ to $SO(V,q)$.  The second has a brief exposition of Clifford algebra and its relation to the pin, spin and orthogonal groups. and some results for these for our $V$ and $q$.   The homomorphism $\psi$ to $Pin(V,q)$ described as a group of units of the Clifford algebra $Cl(V,q)$ is given in the third subsection.   We also show that $\psi$ is a trivial lifting of $\phi$.  Finally, in the last section, we describe a submodule of the Clifford algebra on which the image of $\tilde{y}(m)$ acts in the same way as described by Sidki's $\eta$.

We would like to thank Said Sidki for several productive discussions over the period of his visit to Birmingham in spring 2014 and note that a continuation of this paper will be joint with him.

\section{The groups $y(m)$, $\tilde{y}(m)$ and $Y(m)$}\label{sec:y(m)}

Recall that
\begin{align*}
Y(m)&= \la a_1, \dots, a_m \mid (a_i^ka_j^k)^2=1,i \neq j, k \in \mathbb{Z}\ra,\\
y(m)&= \la a,S_m \mid [s_1,s_1^{a^k}]=1, k \in \mathbb{Z}, a^{s_i}=a^{-1},2\le i \le m-1 \ra
\end{align*}
and $\tilde{y}(m)$ is defined as the semidirect product of $y(m)$ with the cyclic group $\langle \tau \rangle$, where $\tau$ is the involution which inverts $a$ and centralises the $s_i$, for all $i$.  We define $\tilde{s}_i$ to be $\tau s_i$; we note that the $\tilde{s}_i$ also generate a copy $\tilde{S}_m$ of the symmetric group inside $\tilde{y}(m)$.  Let $h(m)$ be the subgroup of $\tilde{y}(m)$ generated by $\tilde{S}_m$ and $\tau$; it is also generated by $S_m$ and $\tau$.

Consider the homomorphism $\pi$ from $\tilde{y}(m)$ onto the quotient obtained by adding the extra relation $a=1$.  From the above presentation, it is easy to see that $\mathrm{Im}(\pi)$ is isomorphic to $S_m \times C_2$ generated by the images of $\tilde{s}_i$, or $s_i$, and $\tau$.  It follows that $\pi$ induces an isomorphism between $h(m)$ and $\mathrm{Im}(\pi) \cong S_m \times C_2$.  In particular, $h(m)$ is a complement to $\ker(\pi)$.  We also note that $\ker(\pi)$ is the normal closure of $a$ in $\tilde{y}(m)$.  Our next goal is to understand the structure of this kernel and we will eventually show that it is isomorphic to $Y(m)$.

We let $b_1 :=a$ and inductively define $b_i := b_{i-1}^{\tilde{s}_{i-1}}$, for $i = 2,\dots,m$.

\begin{lem}\label{homolem}
We have the following
\begin{enumerate}
\item[\emph{(1)}] The involution $\tau$ inverts every $b_i$.
\item[\emph{(2)}] $\tilde{S}_m$ acts on $\{b_1, \dots, b_m\}$ naturally by permuting the indices.
\end{enumerate}
\end{lem}
\begin{pf}
The subgroup generated by $\tilde{s}_2, \dots, \tilde{s}_{m-1}$ centralises $b_1=a$.  This, together with our definitions, implies that the elements $b_1, \dots, b_m$ constitute the orbit under the natural action of $\tilde{S}_m$.  Since $\tau$ commutes with all $\tilde{s}_i$ and inverts $a = b_1$, we see that it inverts every $b_i$.
\end{pf}

Note that, in principle, all the $b_i$ could be equal, however, we will see soon that this cannot happen.

We define a mapping $\theta$ from $Y(m)$ to $\tilde{y}(m)$ by sending every $a_i$ to the corresponding $b_i$.

\begin{prop}
The mapping $\theta$ is a homomorphism.
\end{prop}
\begin{pf}
Since the action of $\tilde{S}(m)$ on $\{b_1, \dots , b_m\}$ is $2$-transitive for $m \geq 3$, it suffices to just check $(b_1^k b_2^k)^2=1$.

\begin{align*}
(b_1^k b_2^k)^2 &= (a^k \tilde{s}_1 a^k \tilde{s}_1)^2 \\
&= (a^k s_1 \tau a^k \tau s_1)^2 \\
&= (a^k s_1 a^{-k} s_1)^2 \\
&=[s_1^{a^{-k}},s_1] \\
&=1
\end{align*}
\end{pf}

By Lemma \ref{homolem}, $\mathrm{Im}(\theta)$ is invariant under the action of $h(m) = \tilde{S}_m \times \la \tau \ra$.  Since the product $\mathrm{Im}(\theta)h(m)$ contains all generators of $\tilde{y}(m)$, we see that $\tilde{y}(m)$ is isomorphic to this product and $\mathrm{Im}(\theta)$ is normal in $\tilde{y}(m)$.  As $\mathrm{Im}(\theta)$ is generated by the conjugates $b_i$ of $a$, $\mathrm{Im}(\theta)$ is the normal closure of $a$ in $\tilde{y}(m)$.  That is, $\mathrm{Im}(\theta) = \ker(\pi)$.  Therefore, we have the following:

\begin{cor}
The group $\tilde{y}(m)$ is isomorphic to the semidirect product of $\mathrm{Im}(\theta)$ with $h(m) = \tilde{S}_m \times \la \tau \ra$.
\end{cor}

We define $\tilde{Y}(m)$ as the semidirect product of $Y(m)$ with the group $S_m \times C_2$, where $S_m$ acts by permuting the $a_i$ naturally and the direct factor $C_2$ inverts all $a_i$.  Let $\tilde{s}_1', \dots, \tilde{s}'_{m-1}$ be the standard generators $(i, i+1)$ of this copy of $S_m$ and $\tau'$ be the generator of the $C_2$.

\begin{lem}
There is a homomorphism $\rho$ from $\tilde{y}(m)$ to $\tilde{Y}(m)$ mapping $a$ to $a_1$, $s_i$ to $s_i':=\tilde{s}_i'\tau'$ and $\tau$ to $\tau'$.
\end{lem}
\begin{pf}
Clearly, $s_i'$ generate a copy of $S_m$ and $\tau'$ commutes with all of them and inverts $a_1$.  Also, it is clear that $a_1^{s_i'} = a_1^{-1}$, for $2 \leq i \leq m-1$.  Hence we just need to see that the commutation relation holds.
\begin{align*}
[\tilde{s}_1'\tau', (\tilde{s}_1'\tau')^{a_1^k}] &= (\tilde{s}_1'\tau' a_1^{-k} \tilde{s}_1'\tau' a_1^k)^2 \\
&= (\tilde{s}_1' a_1^k \tilde{s}_1' a_1^k)^2 \\
&= (a_2^k  a_1^k)^2 \\
&=1
\end{align*}
\qed
\end{pf}
Since $a_i = \rho(a^{\tilde{s}_1 \dots \tilde{s}_{i-1}}) = \rho(b_i)$ for $2 \leq i \leq m$, $\rho$ is surjective.  Moreover, after restricting $\rho$ to $\mathrm{Im}(\pi)$, $\rho$ and $\theta$ are inverses.  In particular, $\theta$ is injective and hence we have the following:

\begin{prop}
The group $\tilde{y}(m)$ is isomorphic to the semidirect product $Y(m):(S_m \times C_2)$.
\qed
\end{prop}

\section{Homomorphisms from $\tilde{y}(m)$}\label{sec:homo}

We will give two homomorphisms from $\tilde{y}(m)$: one into $SO(V,q) = GO(V,q)$ and a second which is the lifting into $Pin(V,q)$, where $V$ is a free module endowed with a quadratic form $q$.  Although the two embeddings are related, the details are sufficiently different and interesting for us to show both.  Before we can do this, we need to define $V$ and $q$.  We begin by letting $R:=\mathbb{F}_2[t,t^{-1}]$, although we will extend this ring later where needed.  Let $V$ be a free module over $R$ of rank $m+1$ with basis $u,v_1,\dots,v_m$.  Our intention is that $\tilde{y}(m)$ acts on $V$ in such a way that the subgroup $S_m$ acts naturally on $v_1,\dots,v_m$ and fixes $u$.  Under this assumption, it can be shown that the form $q$ is uniquely defined up to a scalar.  Namely, we define the symplectic bilinear form $(\cdot,\cdot)$ as follows:
\[
(u,v_i)=(v_i,v_j)=1, \qquad \mbox{for all } i \mbox{ and all } j \neq i.
\]
Associated with this symplectic form, we have a (pseudo-)quadratic form $q$ defined on the basis as follows:
\begin{align*}
q(u) &=1, \\
q(v_1)&= \dots = q(v_m) = t^{-1}.
\end{align*}
This data fully specifies the form $q$ and allows us to compute the value for any vector.  For example, $q(v_i+v_j) = q(v_i) + q(v_j) + (v_i,v_j) = t^{-1} + t^{-1} + 1 = 1$, for all $i \neq j$.

Recall that a \emph{hyperbolic line} is a rank $2$ submodule $W$ spanned by two vectors $e$ and $f$ such that $q(e)=q(f)=0$ and $(e,f)=1$.  For vector spaces $V$ defined over a field, Witt's lemma decomposes $V$ as an orthogonal sum of hyperbolic lines and a rank $1$, or $2$ anisotropic subspace (that is, a subspace which contains no singular vectors).  For modules defined over rings however, no such decomposition exists in general.  However, in our case, we do have a decomposition.

\begin{prop}\label{Vdecomp}
The module $V$ over $\mathbb{F}_2[t,t^{-1}]$ can be decomposed as an orthogonal decomposition $V_1 \perp \dots \perp V_k \perp U$, where $V_i$ are hyperbolic lines and $U$ is a rank $2$ or $3$ submodule.
\end{prop}

\begin{pf}
We will show this by giving an explicit algorithm for the inductive step.  Let $W$ be a free module with basis $e_0,e_1,\dots,e_k$, where $(e_i,e_j)=1$, for all $i\neq j$.  Defining $q_i:=q(e_i)$, we will assume that $q_2=\dots=q_k$.  We further suppose that $W$ has rank at least $4$.  To find a hyperbolic pair $e,f$, we look for solutions $(\alpha,\beta)$ to the equation
\begin{equation}\label{eqn}
0=q(\alpha e_0+\beta e_1 + e_k) = \alpha^2 q_0 + \beta^2 q_1 + q_2 + \alpha + \beta + \alpha \beta.
\end{equation}
If such a solution exists, then set $e= \alpha e_0+\beta e_1 + e_k$ and $f = \alpha e_0+\beta e_1 + e_{k-1}$.  Now, both $e$ and $f$ are singular vectors and
\begin{align*}
(e,f) &= (\alpha e_0+\beta e_1 + e_k, \alpha e_0+\beta e_1 + e_{k-1}) \\
&= \alpha \beta + \alpha + \beta \alpha + \beta + \alpha + \beta + 1\\
&= 1,
\end{align*}
hence they form a hyperbolic pair.

We note that $(e_0 ,e) = (e_0,f) = \beta+1$, $(e_1 ,e) = (e_1,f) = \alpha+1$ and $(e_i ,e) = (e_i,f) = \alpha + \beta+1$, for $i=2, \dots, k-2$.  In particular, $(e_i, e+f)=0$ for $i=0, \dots, k-2$.  Hence, we define $e_0' = e_0 + (\beta+1)(e+f)$, $e_1' = e_1 + (\alpha+1)(e+f)$ and $e_i' = e_i + (\alpha + \beta+1)(e+f)$, for $i = 2, \dots, k-2$.  It is clear that $e_0',e_1', \dots, e_{k-2}'$ span the perp of $\langle e,f \rangle$.  Moreover,
\begin{align*}
(e_0',e_1') &= (e_0 + (\beta+1)(e+f),e_1 + (\alpha+1)(e+f))\\
&=1+0+0+0 \\
&= 1
\end{align*}
and similarly for other $i$, so we have $(e_i',e_j')=1$ for all $i \neq j$.

Hence, provided we can always find a solution to (\ref{eqn}), the above construction gives an inductive step for the proof.  However, the norms $q_i$ of the vectors change in the inductive step.  That is, $q_0':=q(e_0') = q(e_0 + (\beta+1)(e+f)) = q_0 + (\beta+1)^2 q(e+f) = q_0 + \beta^2+1$.  Similarly we have that $q_1':=q(e_1') = q_1+\alpha^2+1$ and $q_i':=q(e_i') = q_i + \alpha^2 + \beta^2 + 1$, for $i = 2, \dots, k-2$.  So, we must show that there is always a solution $(\alpha,\beta)$ to (\ref{eqn}) given the particular values of the $q_i$ at each step.  Starting with $q_0 = 1, q_1 = q_2 = t^{-1}$, we compute the next four solutions and values for $q_i$.

\begin{table}[!ht]
\centering
\begin{tabular}{|c|c|c|c|c|}
\hline
$q_0$ & $q_1$ & $q_2$ & $\alpha$ & $\beta$\\
\hline
$1$&$t^{-1}$&$t^{-1}$&$1$&$1$\\
$1$&$t^{-1}$&$t^{-1}+1$&$0$&$1$\\
$1$&$t^{-1}+1$&$t^{-1}+1$&$1$&$1$\\
$1$&$t^{-1}+1$&$t^{-1}$&$0$&$1$\\
$1$&$t^{-1}$&$t^{-1}$&&\\
\hline
\end{tabular}
\end{table}

Since this has four-fold periodicity and by assumption $W$ has rank at least $4$, this completes the induction to show that $V_1 \perp \dots \perp V_k \perp U$ is a rank $m+1$ submodule and $U$ has rank $2$ or $3$.  Note that, the new basis we produce is a $\mathbb{F}_2$-linear combination of the original basis, hence $V = V_1 \perp \dots \perp V_k \perp U$.
\end{pf}

Note that in the above proposition, we do not claim that $U$ cannot be decomposed further.  In fact, one can show that if $U$ has rank $2$ then it is not a hyperbolic line and we suspect that a detailed case analysis of the rank $3$ case will show that that too does not decompose further.  However, the above is enough for our needs.  Our result above actually decomposes $V$ as an $\mathbb{F}_2$-space; this suggests that we may be able to analyse the module in a different way to produce a better such result.

We further note that when $U$ has rank $3$, it need not be anisotropic.  For example, if $V$ is itself rank $3$, $u+v_1+v_2$ is singular and is the radical of the form and the only singular vector up to multiplication by elements of the ring.  So, $V$ has no hyperbolic lines, but is not anisotropic.

Proposition \ref{Vdecomp} will be useful in determining the structure of the Clifford algebra and the pin and spin groups.  However, we begin by defining a homomorphism into $SO(V,q)$.

\subsection{Homomorphism into the special orthogonal group}\label{orthog}

We need to define the action of $\tau$, $a$ and $\tilde{s}_i := \tau s_i$ on $V$.  Let $w$ be a vector in $V$ such that $q(w) \neq 0$.  Then $r_w$ will denote the transvection given by $r_w(x):=x+\frac{(x,w)}{q(w)}w$.  This is an orthogonal transformation; that is, it preserves $q$ and $(\cdot,\cdot)$.

\begin{prop}\label{orthoghom}
There is a homomorphism $\phi:\tilde{y}(m) \to SO(V,q)$ defined as follows
\begin{align*}
\tau &\mapsto r_u,\\
a &\mapsto r_u r_{v_1}, \\
\tilde{s}_i &\mapsto r_{v_i + v_{i+1}}.
\end{align*}
\end{prop}

As noted before, $q(v_i+v_j) = 1$, for all $i \neq j$, and so our formulae make sense.

The proof will consist of a series of lemmas.

\begin{lem}
The orthogonal transformations $\phi(\tilde{s}_i)$ generate a group $S_m$ which naturally permutes $v_1, \dots, v_m$ and fixes $u$.
\end{lem}
\begin{pf}
Let $w$ equal $u$, or $v_j$ with $j \neq i,i+1$.  Then, $(w,v_i+v_{i+1})=0$, so $r_{v_i+v_{i+1}}(w)=w$.  On the other hand, $r_{v_i+v_{i+1}}(v_i)=v_i+\frac{1}{1}(v_i + v_{i+1}) = v_{i+1}$.  Similarly, for $v_{i+1}$.  This shows that $\phi(\tilde{s}_i)$ swaps $v_i$ and $v_{i+1}$ and fixes all other elements of the basis.
\end{pf}

\begin{lem}
We have that
\begin{enumerate}
\item[\emph{(1)}] $\phi(\tau)$ commutes with $\phi(\tilde{s}_i)$, for all $1 \leq i \leq m-1$,
\item[\emph{(2)}] $\phi(\tilde{s}_i)$ commutes with $\phi(a)$, for all $2 \leq i \leq m-1$,
\item[\emph{(3)}] $\phi(\tau)$ inverts $\phi(a)$.
\end{enumerate}
\end{lem}
\begin{pf}
Since $\phi(\tilde{s}_i)$ fixes $u$, it commutes with $r_u = \phi(\tau)$, giving the first part.  Similarly, $\phi(\tilde{s}_i)$, for $2 \leq i \leq m-1$, fixes $v_1$ and $u$ and hence $\phi(\tilde{s}_i)$ commutes with both $r_u$ and $r_{v_1}$ and hence also with $\phi(a)$.  Since $r_u$ and $r_{v_1}$ both have order two, $r_u$ inverts $\phi(a) = r_ur_{v_1}$.
\end{pf}

It remains to prove the final commutator relation, $[\phi(s_1)^{\phi(a)^k},\phi(s_1)]=1$ for all $k \geq 1$.  We do this via two lemmas:

First, we temporarily extend our ring which will allow us to write $\phi(s_1)^{\phi(a)^k}$ in a much nicer form.  Let $s = \sqrt{t}$ and $\alpha$ and $\alpha^{-1}$ be the roots of $x^2 + sx +1$.  Then $s = \alpha + \alpha^{-1}$ and $t = \alpha^2 + \alpha^{-2}$.
\begin{lem}
The action of $\phi(s_1)^{\phi(a)^k}$ on $V$ is given by the matrix
\[
\left (
\begin{array}{ccc|c}
\alpha^{-2k} + 1 + \alpha^{2k} & \alpha^{-2k} + \alpha^{2k} & \alpha^{-2k} + \alpha^{2k} & \\
\Sigma^{k-1} & \Sigma^{k-1} + 1 & \Sigma^{k-1} & 0\\
\Sigma^k & \Sigma^k & \Sigma^k + 1 & \\
\hline

\begin{array}{c} \rule{0pt}{3ex}\alpha^{-2k} + 1 + \alpha^{2k} \\ \vdots \end{array} & \begin{array}{c} \rule{0pt}{3ex}\Sigma^{k-1} + 1  \\ \vdots \end{array} & \begin{array}{c} \rule{0pt}{3ex}\Sigma^k + 1  \\ \vdots \end{array} & I_{m-2,m-2} 
\end{array}
\right )
\]
for $k \geq 1$, where $\Sigma^k = \Sigma_{i = -k}^k \alpha^{2i}$ and $\Sigma^{-1} :=0$.
\end{lem}
We stress that all the entries in the above matrix are actually in $\mathbb{F}_2[t]$ since it is a product of matrices with entries in $\mathbb{F}_2[t]$.  It is just notationally convenient to use $\alpha$.
\begin{pf}
The proof follows by induction, where both the base case and the inductive step are given by straightforward matrix computations.
\end{pf}
We can now see that both $\phi(s_1)^{\phi(a)^k}$ and $\phi(s_1)$ have the form:
\begin{align*}
u & \mapsto u + f_0(u+v_1+v_2), \\
v_1 & \mapsto v_1 + f_1(u +v_1+v_2), \\
v_2 & \mapsto v_2 + f_2(u+v_1+v_2), \\
v_i & \mapsto v_i + (f_0+1)u + (f_1+1)v_1 + (f_2+1)v_2,
\end{align*}
where $f_0,f_1,f_2$ are coefficients in $R$ with $f_0+f_1+f_2 = 0$, and $3 \leq i \leq m-1$.
\begin{lem}
Any two functions of the above type commute.
\end{lem}
\begin{pf}
Suppose we have two such functions $F$ and $G$ with coefficient functions $f_j$ and $g_j$, respectively.  Note that as $f_0+f_1+f_2 = 0$, $F$ and $G$ fix the vector $w:=u +v_1+v_2$.  They also both fix $v_3+v_4, \dots, v_{m-2}+v_{m-1}$.  Since we can complete these to a basis by adding $v_1$, $v_2$ and $v_3$, we need just consider the action on $v_1$, $v_2$ and $v_3$.  Both $F$ and $G$ act on $v_1$ and $v_2$ by adding a scalar multiple of $w$.  Since they also fix $w$, $F$ and $G$ commute on $v_1$ and $v_2$.  Finally, for $v_3$ we compute $v_3FG$:
\begin{align*}
v_3 &\overset{F}{\mapsto} v_3 + w + f_0u + f_1v_1 + f_2v_2 \\
&\overset{G}{\mapsto} v_3 + w + g_0(u+f_0w) + g_1(v_1+f_1w) + g_2(v_2+f_2w)\\
& \qquad + w + f_0(u+g_0w) + f_1(v_1+g_1w) + f_2(v_2+g_2w)\\
&= v_3+ (f_0+g_0)u +(f_1+g_1)v_1+ (f_2+g_2)v_2
\end{align*}
Hence, $F$ and $G$ commute.
\end{pf}
\begin{cor}
We have that $[\phi(s_1)^{\phi(a)^k},\phi(s_1)]=1$ for all $k \geq 1$.
\end{cor}
This completes the proof of Proposition \ref{orthoghom}.  It is clear from the definition of $\phi$ that $\mathrm{Im}(\phi)$ is contained in the subgroup $TO(V,q) \leq SO(V,q)$ generated by all the orthogonal transvections.

\medskip

\subsection{The Clifford algebra, pin, spin and orthogonal groups}\label{sec:clifford}

Before giving our second homomorphism, we give a brief exposition of the Clifford algebra and the spin, pin and orthogonal groups in general.  We also give some specific results on these for our $V$ and $q$.  For a more detailed discussion, we refer the reader to \cite[Chapter 7]{omeara}.

Let $Cl = Cl(V,q)$ be the Clifford algebra on $V$; that is, the quotient algebra $T(V)/I$ of the tensor algebra $T(V)$ by the ideal $I$ generated by the relations
\[
w^2 = q(w)
\]
for all $w \in V$.  Other useful relations derived from the above are:
\[
ww' + w'w = (w,w')
\]
for all $w,w' \in V$.  Since $V$ and $R$ embed naturally in the algebra $Cl$, we often abuse notation and say that $R$ and $V$ lie in $Cl$.  

Before defining the spin and orthogonal groups, we first note some features of the algebra.  Since the tensor algebra has an $\mathbb{N}$-grading given by the rank of tensors, the quotient $Cl$ inherits a $\mathbb{Z}_2$-grading.  That is, $c \in Cl$ is in the even part, notated by $Cl^0$, if $c$ is the sum of tensors of even rank and in the odd part, $Cl^1$, if it is the sum of tensors of odd rank.

The algebra also has some natural automorphisms.  The map $-1_V:V \to V$ can be extended to an automorphism $\alpha$ on $Cl$.  Note that it is the identity on the even part $Cl^0$ and acts by negation on the odd part (since we work in characteristic two, this map will not concern us: it is included here for the sake of completeness).  There is also a \emph{transpose} map from $Cl$ to $Cl^{op}$ given by reversing the order of vectors in products, i.e., it is defined by extending linearly the map $(w_1 \dots w_k)^{tr} = w_k \dots w_1$, where $w_i \in V$. Clearly, the transpose map is an anti-automorphism of $Cl$.  Combining these two maps, we get the \emph{Clifford conjugation} map $\overline{x} = \alpha(x^{tr})$ which is also an anti-automorphism of $Cl$.  Note that in characteristic two, Clifford conjugation is just the transpose map.

We may now define the Clifford group $C(V,q):= \{ c \in Cl(V,q)^\times : c^{-1} v \alpha(c) \in V \, \forall v \in V\}$, where $Cl(V,q)^\times$ denotes the set of units in $Cl$.  From this definition, there is a natural homomorphism $\pi:C(V,q) \to GO(V,q)$ given by mapping $c$ to the map $v \mapsto c^{-1}v\alpha(c)$.  If $w \in V$ is an anisotropic vector, then it is an invertible element of $Cl$ and $w^{-1} = w/q(w)$.  Since $vw+wv = (v,w)$, we have that
\[
w^{-1}v\alpha(w) = w^{-1}(wv - (v,w)) = v - \frac{(v,w)}{q(w)}w.
\]
Hence, $\pi(w)$ is an orthogonal transvection $r_w$ in $TO(V,q) \leq SO(V,q)$ (without the use of $\alpha$ in the definition it would be $-r_w$).

We can now define $Pin(V,q):=\{ c \in C(V,q) : c \bar{c}=1 \}$ and $Spin(V,q)$ as the even part of $Pin(V,q)$.  The map $x \mapsto x \overline{x}$ is sometimes called the \emph{spinor norm}.  We define $O'(V,q)$ to be the image of $Spin(V,q)$ under the natural map $\pi$; this is sometimes called the \emph{spinorial kernel}.

\medskip

In order to define the second homomorphism, into $Pin(V,q)$, we first need to extend our ring $R$ by adding $s:=\sqrt{t}$ and hence also its inverse $s^{-1}=st^{-1}$; by an abuse of notation we also call this larger ring $R$.  Hence, we now work over the ring $R=\mathbb{F}_2[s,s^{-1}] = \mathbb{F}_2[\sqrt{t},\sqrt{t^{-1}}]$.  The quadratic and associated bilinear forms have the same values on the basis vectors as before and are extended linearly to the larger module.  We use the same notation $V$ for it; similarly $q$ and $(\cdot,\cdot)$.

Before we define our homomorphism, we first briefly discuss the structure of our Clifford algebra $Cl$ and the pin, spin and orthogonal groups.  We begin by quote two lemmas and apply them to our case.

\begin{lem}\emph{\cite[Lemma 7.1.9]{omeara}}
Suppose $V$ has an orthogonal splitting $V=V_1 \perp \dots \perp V_k \perp U$, where the $V_i$ are hyperbolic lines and $U$ is a free submodule of finite rank.  Then,
\[
Cl(V,q) \cong {\rm Mat}_{2^k}(Cl(U,q)),
\]
the $2^k \times 2^k$ matrix algebra with entries in $Cl(U,q)$.
\end{lem}

Hence, using Proposition \ref{Vdecomp}, we see that $Cl(V,q)$ is isomorphic to the matrix algebra with entries in a Clifford algebra $Cl(U,q)$, where $U$ is the rank $2$, or $3$ submodule of $V$ described.

\begin{lem}\emph{\cite[Theorem 7.1.14]{omeara}}\label{Clcentre}
Let $V$ be a free module of finite rank over a ring $R$ with a quadratic form and where the associated bilinear form is non-singular.  Suppose that $V$ has an orthogonal decomposition $V_1 \perp \dots \perp V_k$ into submodules $V_i$ of rank $1$, or $2$.  Then, if $V$ has even rank, the centre $Cen(Cl) = R$ and $Cen(Cl^0)$ is a free module of rank two.  If $V$ has odd rank, then the centre $Cen(Cl)$  is a free module of rank two and $Cen(Cl^0) = R$. 
\end{lem}

Our bilinear form is non-singular precisely when $m+1$ is even, that is when $m$ is odd.  In this case, using Proposition \ref{Vdecomp}, we see that $V$ has the required orthogonal decomposition into submodules of rank $2$.  Hence, when $m$ is odd, the centre of $Cl$ is $R$.

\begin{lem}\label{moddPiniso}
When $m$ is odd, the natural map $\pi:C(V,q) \to SO(V,q)$ restricted to $Pin(V,q)$ and hence also $Spin(V,q)$, is injective.
\end{lem}
\begin{pf}
The kernel of $\pi$ is the set of all $c \in Cl$ which commutes with every $v \in V$.  However, if $c$ commutes with all $v \in V$, then it commutes with all products of vectors and hence lies in the centre of $Cl$.  Since $m$ is odd, this is just $R$.  However, $r\bar{r} = r^2$, hence $R \cap Pin(V,q) = 1$ and $\pi$ is injective.
\end{pf}

When $m+1$ is odd, that is when $m$ is even, the bilinear form on $V$ necessarily has a radical.  In fact, this radical is a rank $1$ submodule spanned by $r:=u+v_1+ \dots + v_m$.  Let $N$ now be the submodule spanned by $u,v_1, \dots, v_{m-1}$, so $V = \la r \ra \perp N$.  Hence, every element $z$ of $Cl(V,q)$ can be written as $x + yr$, where $x,y \in Cl(N,q)$.  If $z \in Cen(Cl(V,q))$, then $z$ commutes with elements of $Cl(N,q)$, hence $x, y \in Cen(Cl(N,q))$.  However, $N$ has even rank, so its centre is just $R$.  We have shown the following:

\begin{lem}
When $m$ is even, $Cen(Cl)$ is a rank $2$ submodule spanned by $1$ and $r=u+v_1+\dots + v_m$.
\end{lem}

\begin{cor}\label{mevenPiniso}
When $m$ is even, $\pi$ restricted to $Pin(V,q)$ is not injective.  Moreover, elements of $\ker \pi \cap Pin(V,q)$ have the form $1+\alpha r$ when $m=2 \mod 4$ and $1+ \alpha (1+r)$ when $m=0 \mod 4$, where $\alpha \in R$.  When restricted to $Spin(V,q)$, however, $\pi$ is injective.
\end{cor}
\begin{pf}
Letting $z = \alpha 1 + \beta r$, we just need to check when $1=z\bar{z} = \alpha^2+\beta^2r^2$.  Since $r^2=q(r)$ and $q(r)$ is 0 when $m=2 \mod 4$ and 1 when $m=0 \mod 4$, the result follows.  Since the only even element in $\ker \pi$ is $1$, $\pi$ is injective when restricted to $Spin(V,q)$.
\end{pf}

Having identified the kernel of $\pi$ when restricted to the pin and spin groups, we now turn our attention to the size of the image in $SO(V,q)$.  For the proofs, we will refer heavily to results from \cite{omeara}.  Let $C^+(V,q)$ denote the even part of $C(V,q)$.

\begin{prop}\label{C^+exact}
The following sequence is exact
\[
1 \to R^* \hookrightarrow C^+(V,q) \stackrel{\pi}{\longrightarrow} SO(V,q) \stackrel{\mathrm{R}}{\longrightarrow} \mathbb{Z}_2 \to 1
\]
\end{prop}
\begin{pf}
By \cite[Theorem 7.2.18]{omeara}, the following sequence is exact
\[
1 \to R^* \hookrightarrow C^+(V,q) \stackrel{\pi}{\longrightarrow} GO(V,q) \stackrel{(\Psi,\mathrm{R})}{\longrightarrow} \mathrm{Pic}_2(R) \times \mathbb{Z}_2(R)
\]
where $\mathrm{Pic}_2(R)$ is the subgroup of the Picard group containing all involutions and $\mathbb{Z}_2(R)$ is the group of idempotents of $R$ (the operation is given by $e_1 \dot{+} e_2 = e_1+e_2-2e_1e_2$).  Since $R=\mathbb{F}_2[s,s^{-1}]$ is a principal ideal domain, the Picard group, and hence $\mathrm{Pic}_2(R)$, is trivial.  Also, it is clear that the only idempotents in $R=\mathbb{F}_2[s,s^{-1}]$ are $0$ and $1$, hence $\mathbb{Z}_2(R) = \mathbb{Z}_2$ for us.  Hence, it remains to show that the map $\mathrm{R}:SO(V,q) \to \mathbb{Z}_2$ is surjective.  Pick $w \in V$ anisotropic and let $r_w$ be the orthogonal transvection with respect to $w$.  By \cite[Example 1, p.241]{omeara}, $\mathrm{R}(r_w)=1$, hence $\mathrm{R}$ is indeed onto.
\end{pf}

Let $O^+(V,q)$ be the kernel of $\mathrm{R}:SO(V,q) \to \mathbb{Z}_2$.  The above result shows that $\pi(C^+(V,q)) = O^+(V,q)$ has index $2$ in $SO(V,q)$.

\begin{prop}\label{Spinexact}
When $m \geq 4$, the following sequence is exact
\[
1 \to Spin(V,q) \stackrel{\pi}{\longrightarrow} O^+(V,q) \stackrel{\Theta}{\longrightarrow} \mathbb{Z}_2 \to 1
\]
\end{prop}
\begin{pf}
By \cite[Theorem 7.2.21]{omeara}, the following sequence is exact
\[
1 \to \mu(R) \hookrightarrow Spin(V,q) \stackrel{\pi}{\longrightarrow} O^+(V,q) \stackrel{\Theta}{\longrightarrow} \mathrm{Disc}(R)
\]
where $\mu(R) = \{r \in R^* | r^2=1\}$ and $\mathrm{Disc}(R)$ is the group of isomorphism classes of discriminant modules of the ring.  Clearly, the only element in $R=\mathbb{F}_2[s,s^{-1}]$ which squares to $1$ is $1$ itself.  For  $\mathrm{Disc}(R)$, we have the following exact sequence:
\[
1 \to (R^*)^2 \hookrightarrow R^* \to \mathrm{Disc}(R) \to \mathrm{Pic}_2(R) \to 1
\]
So, since $\mathrm{Pic}_2(R) = 1$, $\textrm{Disc}(R) \cong {R^*}/{(R^*)^2}$ which is isomorphic to $\mathbb{Z}_2$ for $R=\mathbb{F}_2[s,s^{-1}]$.  It remains to show that $\Theta$ is surjective.  Since $m \geq 4$, by Proposition \ref{Vdecomp}, $V$ contains a hyperbolic line $(e,f)$.  Let $w = sef+fe$.  By an easy calculation, $w^{-1} = s^{-1}ef + fe$.  So, $w^{-1}ew = s^{-1}e$, $w^{-1}fw = sf$ and $w^{-1}xw = x$ for $x \in \{e,f\}^\perp$.  Since $w$ is clearly even, $w \in C^+(V,q)$.  By \cite[Example 3, p. 241]{omeara}, $\mathrm{R}(\pi(w)) = 0$, hence $\pi(w) \in O^+(V,q)$.  Furthermore, $\Theta(\pi(w)) = w\bar{w}(R^*)^2$.  Since $w\bar{w} = (sef + fe)(sfe + ef) = s(ef+fe) = s$ and $s \not \in (R^*)^2$, $\Theta$ is onto.
\end{pf}

\begin{cor}
The \emph{(}image of\emph{)} $Spin(V,q)$ has index $4$ in $SO(V,q)$.
\end{cor}
\begin{pf}
By Lemma \ref{moddPiniso} and Corollary \ref{mevenPiniso}, $\pi$ is injective when restricted to $Spin(V,q)$.  Hence, $O'(V,q):=\pi(Spin(V,q)) \cong Spin(V,q)$.  However, by Propositions \ref{C^+exact} and \ref{Spinexact}, $O'(V,q)$ has index $2$ in $O^+(V,q)$ which in turn has index $2$ in $SO(V,q)$.
\end{pf}

\subsection{Homomorphism into the Clifford group}\label{sec:pinhom}

We now describe our homomorphism.  Recall that we have extended our ring $R$ to $R=\mathbb{F}_2[s,s^{-1}]$.  Using the relations in the Clifford algebra, we have the following equalities for our basis vectors
\begin{align*}
u^2 &=1, \\
v_i^2 &= t^{-1}, \\
(v_i+v_j)^2 &=1,\\
uv_i + v_i u &= 1, \\
v_iv_j + v_jv_i &= 1,
\end{align*}
for all $i \neq j$.  In particular, $u$, $v_1$ and $v_i + v_{i+1}$ are all invertible and $u^{-1} = u$, $v_1^{-1} = t v_1$ and $(v_i + v_{i+1})^{-1} = v_i + v_{i+1}$.  Note that although $v_1$ is not an involution, $sv_1$ is.

\begin{prop}
The map $\psi:\tilde{y}(m) \to Pin(V,q)$ given by
\begin{align*}
\tau &\mapsto u, \\
a &\mapsto s u v_1,\\
\tilde{s}_i &\mapsto v_i + v_{i+1},
\end{align*}
on the generators of $\tilde{y}(m)$ defines a group homomorphism.
\end{prop}

As noted above, $u$ and $v_i + v_{i+1}$ are involutions.  The following easy lemma will be useful in verifying the remaining relations.

\begin{lem}\label{s_i commute}
We have
\begin{enumerate}
\item[\emph{(1)}] $(v_i+v_{i+1})u = u(v_i+v_{i+1})$
\item[\emph{(2)}] $(v_i+v_{i+1})v_j = v_j(v_i+v_{i+1})$ provided $j \neq i,i+1$
\item[\emph{(3)}] $(v_i+v_k)v_i = v_k (v_i+v_k)$
\end{enumerate}
\end{lem}
\begin{pf}
For the first identity, $(v_i+v_{i+1})u = v_iu+v_{i+1}u = 1+uv_i+1+uv_{i+1} = u(v_i+v_{i+1})$.  The remaining calculations are similar.
\end{pf}

We write $w_i :=v_i+v_{i+1} = \psi(\tilde{s}_i)$.   It is clear from the above lemma that $\psi(\tilde{s}_i)$ commutes with $\psi(\tilde{s}_j)$ when $|i-j| \geq 2$.  Observe that from the Clifford algebra, we have $w_iw_{i+1} + w_{i+1}w_i = (w_i,w_{i+1}) = 1$.  So, we get
\begin{align*}
\big(\psi(\tilde{s}_i)\psi(\tilde{s}_{i+1})\big)^3 &:= (w_iw_{i+1})^3\\
&=(1+w_{i+1}w_i)w_iw_{i+1}w_iw_{i+1}\\
&=w_iw_{i+1}w_iw_{i+1} +w_iw_{i+1}\\
&=(1+w_{i+1}w_i)w_iw_{i+1}+w_iw_{i+1}\\
&=1
\end{align*}

This verifies the Coxeter presentation for the group $\tilde{S}_m \le \tilde{y}(m)$ which means that $\psi$ is a well-defined isomorphism when restricted to the subgroup $\la \tilde{s}_1, \dots, \tilde{s}_{m-1} \ra$. It is also clear that $\psi(\tau)$ inverts $\psi(a)$ and centralises the $\psi(\tilde{s}_i)$.  To show the other relations involving $a$, we first need the following lemma.

Using the above relations, when $j \neq 1$, we have
\[
\psi(a)^{\psi(s_j)} = s (v_j + v_{j+1})uuv_1 u(v_j + v_{j+1}) = s v_1u = \psi(a^{-1}).
\]

So, it remains to prove the commutation relation $[\psi(s_1),\psi(s_1)^{\psi(a)^k}]=1$.  We begin by using Lemma \ref{s_i commute} to get
\[
\psi(s_1) \psi(s_1)^{\psi(a)^k} = t^k u(v_1+v_2) (v_1 u)^k u (v_1+v_2) (uv_1)^k = t^k (uv_2)^k(uv_1)^k
\]
and similarly $\psi(s_1)^{\psi(a)^k}\psi(s_1) = t^k (v_1u)^k(v_2u)^k$.  Hence it remains to verify that $(v_1u)^k(v_2u)^k=(uv_2)^k(uv_1)^k$ for all $k\ge 0$. 

\begin{lem} \label{relation}
The algebra element sequences $(uv_i)^k$ and $(v_iu)^k$ satisfy the recurrence  
relation $\alpha_k=\alpha_{k-1}+t^{-1}\alpha_{k-2}$.
\end{lem}

\begin{pf}
Note that $(uv_i)^k=(uv_i)^{k-2}(uv_i)(v_iu+1)=(uv_i)^{k-2}(t^{-1}+uv_i)=
(uv_i)^{k-1}+t^{-1}(uv_i)^{k-2}$.  The same calculation also works for $(v_iu)^k$.
\end{pf}

Let $a_k$ and $b_k$ be the sequences of elements of $R$ satisfying the 
above recurrence relation and the initial conditions $a_0=1$, $a_1=0$ and 
$b_0=0$, $b_1=1$, respectively. These sequences can be computed explicitly, 
but we do not really need those formulas and so we skip the computation.

\begin{lem} \label{expressions}
We have
$$(uv_i)^k=a_k+b_kuv_i=(a_k+b_k)+b_kv_iu$$
and
$$(v_iu)^k=a_k+b_kv_iu=(a_k+b_k)+b_kuv_i.$$
\end{lem}

\begin{pf}
Again, we just show the first claim. Note that $(uv_i)^0=1=a_0+b_0uv_i$ 
and $(uv_i)^1=uv_i=a_1+b_1uv_i$. Now in view of Lemma \ref{relation}, we obtain 
inductively, for $k\ge 2$, that $(uv_i)^k=(uv_i)^{k-1}+t^{-1}(uv_i)^{k-2}=
(a_{k-1}+b_{k-1}uv_i)+t^{-1}(a_{k-2}+b_{k-2}uv_i)=(a_{k-1}+t^{-1}a_{k-2})+
(b_{k-2}+t^{-1}b_{k-2})uv_i=a_k+b_kuv_i$. Clearly, $a_k+b_kuv_i=
(a_k+b_k)+b_kv_iu$, since $v_iu=uv_i+1$.
\end{pf}

We will also need the following observation.

\begin{lem} \label{commutation}
If $i\ne j$ then $(au+bv_i)(cu+dv_j)=(cu+dv_j)(au+bv_i)+(ad+bc+bd)$.
\end{lem}

\begin{pf}
We have $(au+bv_i)(cu+dv_j)=acu^2+aduv_j+bcv_iu+bdv_iv_j=acu^2+ad(v_ju+1)+
bc(uv_i+1)+bd(v_jv_i+1)=(cau^2+dav_ju+cbuv_i+dbv_jv_i)+(ad+bc+bd)=
(cu+dv_j)(au+bv_i)+(ad+bc+bd)$.
\end{pf}

\medskip
Finally, we establish our main claim.

\begin{prop}
We have that $(v_1u)^k(v_2u)^k=(uv_2)^k(uv_1)^k$ for all $k\ge 0$.
\end{prop}

\begin{pf}
By Lemma \ref{expressions}, $(v_1u)^k(v_2u)^k=(a_k+b_kv_1u)((a_k+b_k)+
b_kuv_2)=(a_ku+b_kv_1)u^2((a_k+b_k)u+b_kv_2)=(a_ku+b_kv_1)((a_k+b_k)u+b_kv_2)$. 
In view of Lemma \ref{commutation}, $(a_ku+b_kv_1)((a_k+b_k)u+b_kv_2)=
((a_k+b_k)u+b_kv_2)(a_ku+b_kv_1)+(a_kb_k+b_k(a_k+b_k)+b_k^2)=
((a_k+b_k)u+b_kv_2)(a_ku+b_kv_1)$, since the second summand is obviously zero.

On the other hand, by Lemma \ref{expressions}, we have that 
$((a_k+b_k)u+b_kv_2)(a_ku+b_kv_1)=((a_k+b_k)u+b_kv_2)u^2(a_ku+b_kv_1)=
((a_k+b_k)+b_kv_2u)(a_k+b_kuv_1)=(uv_2)^k(uv_1)^k$, as claimed.
\end{pf}

This shows that $\psi$ is a homomorphism.   To complete the proof we observe that $u\overline{u} = u^2 =1$, $suv_1 \overline{suv_1} = t uv_1 v_1u = 1$ and $(v_i + v_{i+1})\overline{(v_i + v_{i+1})} = (v_i + v_{i+1})^2 =1$, so the image of $\psi$ does indeed lie in $Pin(V,q)$.

We note that when restricted to $y(m)$, all the elements are in the even part of $C$ and so $\psi$ maps $y(m)$ into $Spin(V,q)$.

Recall that $\pi:Pin(V,q) \to SO(V,q)$ is a homomorphism which maps $w$ to the map $v \mapsto w^{-1}v w$.  If $w$ is an anisotropic vector, then $w$ is sent to $r_w$, an orthogonal transvection.

\begin{prop}
The map $\pi$ is injective when restricted to $\psi(\tilde{y}(m))$ and $\phi = \pi \circ \psi$.  Hence, $\psi$ is a trivial lifting of $\phi$.
\end{prop}
\begin{pf}
 It is clear from the definitions of $\phi$ and $\psi$ that $\phi = \pi \circ \psi$.  Consider first $\pi$ restricted to $\psi(y(m))$.  As noted above, $\psi(y(m)) \le Spin(V,q)$ and, by Lemmas \ref{moddPiniso} and \ref{mevenPiniso}, $\pi$ restricted to $Spin(V,q)$ is injective.  Hence, $\pi$ is injective on $\psi(y(m))$.  Since $|\tilde{y}(m):y(m)| = |\psi(\tilde{y}(m)):\psi(y(m))| = |\phi(\tilde{y}(m)):\phi(y(m))| =2$ and $\phi = \pi \circ \psi$, $\pi$ is injective when restricted to $\psi(\tilde{y}(m))$.
\end{pf}
 The above proposition shows that our two homomorphism carry the same information.  However, they do give different actions of our group.

\section{Sidki's action}

Recall that in \cite[Section 4]{si1} Sidki gives a representation of $y(m,n)$.
This can be easily extended to a representation $\eta = \eta_m:y(m) \to SL(2^{m-2},R)$
of $y(m)$, for $m \geq 3$.  When $m=3$, it is given by
\[
a \mapsto \left ( \begin{array}{cc} \alpha & 0 \\0 & \alpha^{-1} \end{array}
\right ) =: a^{(3)}, \quad 
s_1 \mapsto \left ( \begin{array}{cc} 1 & 0 \\1 & 1 \end{array} \right ) =:s_1^{(3)}, \quad 
s_2 \mapsto \left ( \begin{array}{cc} 0 & 1 \\1 & 0 \end{array} \right )=:s_2^{(3)},
\]
where $\alpha$ is an indeterminate.  For $m\geq 4$, we tensor up to obtain
\begin{align*}
a &\mapsto \left ( \begin{array}{cc} a^{(m-1)}&0 \\0 & \big(a^{(m-1)}\big)^{-1} \end{array} \right )=:a^{(m)}, 
& s_1 &\mapsto \left ( \begin{array}{cc} I_{m'}&0 \\I_{m'} & I_{m'} \end{array}
\right )=:s_1^{(m)}, \\ 
s_2 &\mapsto \left ( \begin{array}{cc} 0 & I_{m'} \\I_{m'} & 0\\ \end{array} \right
)=:s_2^{(m)}, & s_3 &\mapsto \left( \begin{array}{cc}
    s_2^{(m-1)}&0\\I_{m'}&s_2^{(m-1)} \end{array} \right) =:s_3^{(m)}, \\
s_k &\mapsto \left( \begin{array}{cc} s_{k-1}^{(m-1)}
    &0\\0&s_{k-1}^{(m-1)} \end{array} \right) =:s_k^{(m)}
\end{align*}
where $4 \leq k \leq m-1$, $m':=2^{m-3}$ and $I_{m'}$ is the $m'\times m'$ identity matrix.

\medskip

In this section, we will describe a submodule of the Clifford algebra on which $y(m)$ gives this action.  Moreover, we will describe a submodule of rank $2^{m-1}$ (twice the size) where we have an action of $\tilde{y}(m)$ which reduces to Sidki's action.  In other words, inside the Clifford algebra we can both find Sidki's action and extend it to include the non-trivial action of $\tau$.

\begin{lem}
The minimum polynomial of $\psi(a)=s uv_1$ is $x^2 + sx + 1$.
\end{lem}
\begin{pf}
We have $(suv_1)^2+s.suv_1 + 1 = t(1+v_1u)uv_1 + tuv_1 + 1 = 0$.
\end{pf}
Let $\alpha$ and $\alpha^{-1}$ be the two roots of the minimum polynomial.  We now extend our ring $\mathbb{F}_2[s,s^{-1}]$ by adding the two roots; similarly to before, we abuse notation by calling this new ring $R$.

Observe that in the action given above, the first basis vector is an $\alpha$-eigenvector for $a$ and a $1$-eigenvector for $s_1$.  Since we wish our construction to be canonical for all $m \geq 3$, we may look for such a vector $w$ in the subalgebra spanned by $u$, $v_1$ and $v_2$.  Furthermore, by possibly multiplying by $u$, we may look for such a vector in the even part of this subalgebra.  That is, in $U:=\langle 1,uv_1,uv_2,v_1v_2 \rangle$.

\begin{lem}
The vector $w:=\frac{\alpha^2}{s}+\alpha uv_1+\alpha^{-1}uv_2 + sv_1v_2$ is the unique vector in $U$ which is an $\alpha$-eigenvector for $a$ and a $1$-eigenvector for $s_1$, up to scalar multiplication.
\end{lem}
\begin{pf}
%
%

Let $z$ be an $\alpha$-eigenvector for $suv_1$.  Suppose it is also a $1$-eigenvector for $s_1$.  Then combining $zsuv_1=\alpha z$ and $z(uv_1+uv_2)=z$, we get $szuv_2 = (sz+zsuv_1) = (sz+\alpha z) = (s+\alpha)z = \alpha^{-1}z$.  So, $z$ is additionally a $1$-eigenvector for $s_1$ if and only if it is a $\alpha^{-1}$-eigenvector for $suv_2$. 

Let $z=\lambda_0 + \lambda_1uv_1 + \lambda_2uv_2+\lambda_3v_1v_2$ be a vector in $U$.  We compute:
\begin{align*}
za&=(\lambda_0 + \lambda_1uv_1 + \lambda_2uv_2+\lambda_3v_1v_2)suv_1 \\
&= \lambda_0suv_1 + \lambda_1suv_1(1+v_1u) + \lambda_2s(1+v_2u)uv_1+\lambda_3s(1+v_2v_1)(1+v_1u) \\
&= \lambda_0suv_1 + \lambda_1suv_1 + \lambda_1s^{-1} + \lambda_2suv_1 + \lambda_2sv_2v_1\\
& \qquad+\lambda_3s(1+v_2v_1+v_1u+t^{-1}v_2u) \\
&= (\lambda_1s^{-1} + \lambda_2 s+ \lambda_3(s+s^{-1})) +(\lambda_0 +\lambda_1+\lambda_2+\lambda_3)suv_1 +\lambda_3s^{-1}uv_2 \\
&\qquad + (\lambda_2+\lambda_3)sv_1v_2
\end{align*}
We require that $z$ is an $\alpha$-eigenvector for $uv_1$.  By equating coefficients we get
\begin{align*}
\lambda_1s^{-1} + \lambda_2s+\lambda_3(s+s^{-1}) &=\alpha\lambda_0 \\
(\lambda_0+\lambda_1+\lambda_2+\lambda_3)s &=\alpha \lambda_1 \\
\lambda_3s^{-1} &=\alpha\lambda_2 \\
(\lambda_2+\lambda_3)s &=\alpha \lambda_3
\end{align*}
We note that the last two equations are equivalent.

We also require that $z$ is an $\alpha^{-1}$-eigenvector for $uv_2$.  After noticing that $z=\lambda_0 + \lambda_1uv_1 + \lambda_2uv_2+\lambda_3v_1v_2 = (\lambda_0 + \lambda_3) + \lambda_1uv_1 + \lambda_2uv_2+\lambda_3v_2v_1$, we may reuse the above calculation, by interchanging the role of $v_1$ and $v_2$.  Equating coefficients again we get another four equations, including $\lambda_3 = \alpha^{-1}s \lambda_1$.  Hence, $\lambda_2=\alpha^{-2}\lambda_1$.  Using the second equation above, we get $\lambda_0 = s^{-1}(\alpha\lambda_1+s\lambda_1+s\lambda_2 + s\lambda_3) = s^{-1}(\alpha + s+ \alpha^{-2}s+\alpha^{-1}t)\lambda_1 = \frac{\alpha}{s}\lambda_1$.  Using $\lambda_1=\alpha$ gives $w$.  One can check that this is consistent with all the other equations, hence this is the unique solution, up to scalar multiplication.
\end{pf}

Now that we have the vector $w$, we can identify a rank $2^{m-2}$ submodule $W$ of the Clifford algebra $C$ on which our representation of $y(m)$ acts in the same way as Sidki describes.  We define an ordered set $X_m$ which will be a basis for $W$ using the following algorithm.  We begin with $w$ and add $ws_{m-1}$ to the list.  Next we add $ws_{m-2}$ and $ws_{m-1}s_{m-2}$ to our list.  At the $k$\textsuperscript{th} stage, we add $xs_{m-k}$ for each $x$ already in our list.  We continue until we last apply $s_2$.  Since at each of the $m-2$ stages we double the number of vectors, this yields an ordered set $X_m$ of elements of $C$ of size $2^{m-2}$.  It consists of all elements of the form $ws$, where $s$ is an ordered product of the elements in an ordered subset of $\{ s_{m-1}, \dots, s_2 \}$.

\begin{prop}\label{Wbasis}
 The ordered set $X_m=\{ w, ws_{m-1}, \dots, ws_{m-1}\dots s_2 \}$ is linearly independent.
\end{prop}

We will prove the above proposition using induction on $m$ and a number of lemmas.

\begin{lem}
When $m=3$, $w$ and $ws_2$ are linearly independent.
\end{lem}
\begin{pf}
If not, then $w$ is an eigenvector for $s_2$ for some eigenvalue $\lambda$.  Recall that $w$ is an $s^{-1}\alpha^{-1}$ eigenvector for $uv_2$.  We have $\lambda w=ws_2=w(uv_2+uv_3) = s^{-1}\alpha^{-1}w + wuv_3$.  Hence, $w$ is also an eigenvector for $uv_3$.  However, since $w \in U$ and $v_3 \not \in U$, the product $wuv_3$ cannot lie in $U$, a contradiction.
\end{pf}

In order to show the inductive step, it will be easier to consider a slightly different set with a different ordering.

\begin{lem}\label{s_i comms}
 We have that
\[
 s_is_j = \left \{ \begin{array}{ll} 1 & \mbox{if } i=j, \\ s_js_i & \mbox{if } |i-j|\geq 2, \\ s_js_i +1 & \mbox{if } |i-j| =1. \end{array} \right.
\]
\end{lem}
\begin{pf}
 This is a straightforward calculation using Lemma \ref{s_i commute}.
\end{pf}

Define $Y_m$ to be the set comprising all elements of the form $ws$, where $s$ is an ordered product of the elements in an ordered subset of $\{ s_2, \dots, s_{m-1} \}$.

\begin{cor}
 The set $X_m$ is linearly independent if and only if $Y_m$ is linearly independent.
\qed
\end{cor}

We now use induction to show that $Y_{m+1}$ is linearly independent.  Observe that, as an unordered set, $Y_{m+1} = Y_m \cup Y_ms_m$.  We give $Y_m$ a new ordering as follows: we order first by the length of the word $s$ in $ws$ and within that lexicographically.  This ordering induces an ordering on $Y_ms_m$ in the natural way.

Starting with $Y'=Y_m$, we add elements $y$ one at a time from the ordered set $Y_ms_m$ to form a new $Y'$.  We claim that at each step $Y' \cup \{y\}$ is linearly independent.  Indeed, suppose $y:=ws_{i_1}\dots s_{i_l}s_m \in Y_ms_m$.  It has in its decomposition the product $u^lv_1v_{i_1+1}\dots v_{i_l+1}v_{m+1}$, for some $l$.  Since $s_{i_j} = u(v_{i_j}+v_{i_j+1})$ and $Y_m$ is ordered lexicographically, this product does not occur in the decomposition of any element of $Y'$.  Hence $Y' \cup \{y\}$ is linearly independent.  Therefore, by induction, $Y_{m+1}$ is also linearly independent and Proposition \ref{Wbasis} is proved.

\medskip

Define $W$ to be the submodule spanned by $X_m$.

\begin{prop}
The action of $\psi(y(m))$ on $W$ with respect to the ordered basis $X_m$ is as given by Sidki.
\end{prop}
\begin{pf}
We use induction,  starting with the base case, $m=3$.  It is clear that $s_2$ has the required action.  Since $s_2$ inverts $a$ and $w$ is an $\alpha$-eigenvector for $a$, the action of $a$ is also clear.  By choice, $ws_1=w$, so it remains to consider $ws_2s_1$.  Using Lemma \ref{s_i commute} we have
\begin{align*}
ws_2s_1 &= w((v_2+v_3)uu(v_1+v_2)) \\
&=w(1+(v_1+v_2)(v_2+v_3)) \\
&=w+w(v_1+v_2)uu(v_2+v_3) \\
&=w+ws_1s_2 \\
&=w+ws_2
\end{align*}
So, we have shown the case $m=3$.

For the inductive step, consider the two bases: $X_m =\{ w, ws_{m-1},\allowbreak ws_{m-2}, \allowbreak ws_{m-1}s_{m-2}, \dots, ws_{m-1}\dots s_2 \}$ and $X_{m+1}=\{ w, \allowbreak ws_m, \allowbreak ws_{m-1}, \allowbreak ws_ms_{m-1}, \dots, \allowbreak ws_m \dots s_3, ws_2, \dots, \allowbreak ws_m\dots s_2 \}$.  We note that the map on the words $s$ which is generated by $s_i \mapsto s_{i+1}$ maps $X_m$ to the first $2^{m-2}$ basis vectors of $X_{m+1}$.  Therefore, for $k \geq 3$, the action of $\eta_{m+1}(s_k)$ on a basis vector $ws_{i_l}s_{i_{l-1}}\dots s_{i_1}$ in the first half of $X_{m+1}$ is the same as the action of $\eta_m(s_{k-1}) = s^{(m)}_{k-1}$ on $ws_{{i_l}-1}s_{i_{l-1}-1}\dots s_{{i_1}-1}$.  That is, the top half of the action matrix for $\eta_{m+1}(s_k)$ is given by $(s^{(m)}_{k-1} | 0)$.  Now, provided $k \geq 4$, $s_k$ commutes with $s_2$.  So the action of $\eta_{m+1}(s_k)$ on the second half of the basis vectors $ws_2, \dots , ws_m\dots s_2$ is the same as on the first half.  Therefore, the action of $\eta_{m+1}(s_k)$ is as given.  If $k=3$, then $s_2s_3 = 1+s_3s_2$ by Lemma \ref{s_i comms}.  Each basis vector in the second half of $X_{m+1}$ is of the form $wss_2$, where $s$ is a word in $\{s_3, \dots, s_m\}$.  Since $wss_2s_3 = ws + ws s_3s_2$, the action of $\eta_{m+1}(s_3)$ is as given.

For $s_1$, we use the same argument that $wss_2s_1 = ws + ws s_1s_2$ and observe that $s_1$ commutes with all $s_3, \dots, s_m$.  The action of $s_2$ is clearly as given.  Finally, for $a$, recall that $s_i$ inverts $a$ for $i \neq 1$; that is, $s_ia = a^{-1}s_i$.  Hence, an element $ws$ is an $\alpha$-eigenvector if $s$ has even length and an $\alpha^{-1}$-eigenvector if it has odd length.  The first $2^{m-2}$ basis vectors in $X_{m+1}$ have the same length as those in $X_m$, whilst the second half have opposite length parity to those in $X_m$.
\end{pf}

Recall that we may extend $y(m)$ to $\tilde{y}(m)$ by adding an element $\tau$ which inverts $a$ and centralises $s_1, \dots, s_{m-1}$.

\begin{cor}
The action of $\psi(\tilde{y}(m))$ on the rank $2^{m-1}$ submodule $W \oplus Wu$ is given by
\begin{align*}
u & \mapsto \left( \begin{array}{cc} 0 & I_{2^{m-2}} \\  I_{2^{m-2}}& 0 \end{array} \right)\\
 a & \mapsto \left( \begin{array}{cc} \eta_m(a) & 0 \\ 0 & \eta_m(a)^{-1} \end{array} \right) \\
s_i & \mapsto \left( \begin{array}{cc} \eta_m(s_i) & 0 \\ 0 & \eta_m(s_i) \end{array} \right) \qquad \textup{for }i=1,\dots,m-1
\end{align*}
\end{cor}
\begin{pf}
By considering different products in a similar way to the proof of Proposition \ref{Wbasis}, we see that $X_m \cup X_mu$ is linearly independent and hence $W \oplus Wu$ has rank $2^{m-1}$.  The action of $\tilde{y}(m)$ is clear.
\end{pf}

\end{document}